\documentclass[12pt]{article}

\usepackage{amsfonts}
\usepackage{amsmath}
\usepackage{amssymb}
\usepackage{graphicx}
\usepackage[usenames]{color}

\usepackage{mathrsfs}

\newtheorem{theorem}{Theorem}[section]
\newtheorem{proposition}[theorem]{Proposition}

\newtheorem{lemma}[theorem]{Lemma}
\newtheorem{remark}[theorem]{Remark}
\newtheorem{definition}[theorem]{Definition}

\numberwithin{equation}{section} 

\def\Sc{{Schr\"odinger} }

\def\en t{{{\rm Z}\mkern-5.5mu{\rm Z}}}

\def\<{\left<}
\def\>{\right>}
\def\({\left(}
\def\){\right)}

\def\9{{\infty}}
\def\barr{\begin{array}}
\def\earr{\end{array}}

\def\wt{\widetilde}
\def\wh{\widehat}

\def\lbb{{\lambda}}

\def\a{{\alpha}}

\def\n{\noindent }

\def\3{\subset }
\def\na{{\nabla}}
\def\sk{\smallskip }

\def\bk{\bigskip }
\def\e{{\epsilon}}

\begin{document}

\begin{center}
{\Large{\bf Stochastic nonlinear \Sc equations: \\
no blow-up in the non-conservative case}}
\bigskip\bk

{\large{\bf Viorel Barbu}}\footnote{Octav Mayer Institute of
Mathematics (Romanian Academy)   and Al.I. Cuza University and,
700506, Ia\c si, Romania. This work was supported by the DFG through
CRC 701.}, {\large{\bf Michael R\"ockner}}\footnote{Fakult\"at f\"ur
Mathematik, Universit\"at Bielefeld,  D-33501 Bielefeld, Germany.
This research was supported by the DFG through CRC 701.},
{\large{\bf Deng Zhang}}\footnote{Department of Mathematics,
Shanghai Jiao Tong University, 200240 Shanghai, China.}
\end{center}

\bk\bk\bk

\begin{quote}
\n{\small{\bf Abstract.} This paper is devoted to the study of noise
effects on blow-up solutions to stochastic nonlinear Schr\"odinger
equations. It is a continuation of our recent work \cite{BRZ14},
where the (local) well-posedness is established in $H^1$, also in
the non-conservative critical case. Here we prove that in the
non-conservative focusing mass-(super)critical case, by adding a
large multiplicative Gaussian noise,  with high probability one can
prevent the blow-up on any given bounded time interval $[0,T]$,
$0<T<\9$. Moreover, in the case of  spatially independent noise, the
explosion even can be prevented with high probability on the whole
time interval $[0,\9)$. The noise effects obtained here are
completely different from those in
the conservative case studied in \cite{BD03}. } \\


\end{quote}

\vfill

\section{Introduction and main results.}

We consider the stochastic nonlinear Schr\"{o}dinger equation with
linear  multiplicative noise,
\begin{align} \label{equax}
      &idX(t,\xi)=\Delta X(t,\xi)dt+\lambda|X(t,\xi)|^{\alpha-1}X(t,\xi)dt \nonumber \\
      &\qquad \qquad \quad -i\mu(\xi) X(t,\xi)dt+iX(t,\xi)dW(t,\xi),\
      t\in(0,T),\ \xi\in \mathbb{R}^d,
      \\
      &X(0)=x \in H^1.  \nonumber
\end{align}
Here, the exponents of particular interest lie in the focusing
mass-(super)critical range, namely,
\begin{align} \label{f-super}
    \lambda=1,\ \ \a\in[1+\frac{4}{d}, 1+\frac{4}{(d-2)^+}).
\end{align}
$W$ is the colored Wiener process
\begin{equation} \label{W}
     W(t,\xi)=\sum^N\limits_{j=1}\mu_j e_j(\xi)\beta_j(t),\ t\geq
     0,\ \xi\in \mathbb{R}^d,
\end{equation}
where $N<\9$, $\mu_j \in \mathbb{C}$, $e_j$ are real-valued
functions, and $ \beta_j(t)$ are independent real Brownian motions
on a probability space $(\Omega,\mathcal {F},\mathbb{P})$ with
natural filtation $(\mathcal {F}_t)_{t\geq 0}$, $1\leq j \leq N$.
Moreover, as required by the physical context (see \cite{BG09} and
\cite{BPP10}), $\mu$ is of the form
\begin{align}\label{mu}
\mu=\sum\limits_{j=1}^N|\mu_j|^2e_j^2.
\end{align}
Hence $|X(t)|^2_2$ is a martingale, which allows to define the
so-called ''physical probability law''. In particular, in the
conservative case (i.e. $Re \mu_j=0$, $1\leq j\leq N$), the last two
terms in \eqref{equax} coincide with the Stratonovitch integration.
We also refer to \cite{BRZ13} for discussions on the physical
background.

\begin{definition} \label{def-x}
A solution $X$ to \eqref{equax} on $[0,\tau]$, where $\tau$ is an
$(\mathscr{F}_t)$-stopping time, is an $H^1$-valued continuous
$(\mathscr{F}_t)$-adapted process, such that $|X|^{\alpha-1}X\in
L^1(0,\tau;H^{-1})$, $\mathbb{P}-a.s$, and it satisfies
$\mathbb{P}-a.s$
\begin{align} \label{weakx}
    X(t)
    =&x-\int_0^t(i\Delta X(s)+\mu X(s)+\lambda
    i|X(s)|^{\alpha-1}X(s))ds \nonumber \\
     &  +\int_0^t X(s)dW(s),\ \ t\in [0,\tau],
\end{align}
as an equation in $H^{-1}$.
\end{definition}

The well-posedness of \eqref{equax} is studied in our recent paper
\cite{BRZ14}, based on the rescaling transformation used in
\cite{BRZ13} and the Strichartz estimates established in
\cite{MMT08} for perturbations of the Laplacian. We also refer to
the standard monographs \cite{C03} and \cite{LP09} for the
deterministic case (i.e. $\mu_j=0$, $1\leq j\leq N$) and to
\cite{BD03} and \cite{BF12} for the stochastic conservative case
(i.e. $Re \mu_j=0$, $1\leq j\leq N$).

The main interest of this article is to study the noise effects on
blow-up in the focusing mass-(super)critical case. Our motivations
mainly come from two aspects. On the one hand, the blow-up
phenomenon in the deterministic case is extensively studied in the
literature, and it is well known that there exist blow-up solutions
in the focusing mass-(super)critical case \eqref{f-super},
especially for initial data with negative Hamiltonian (cf. e.g.
\cite{C03}, \cite{LP09}). On the other hand, when there is noise in
the system, it is of great interest to investigate the noise effects
on the formation of singularities. For example, in the conservative
case, it is proved in \cite{BD05} in the supercritical case that
noise can accelerate blow-up with positive probability. But in the
critical case numerical results suggest that noise has the effect to
delay explosion (cf. \cite{BDM02}, \cite{DM02.1} and \cite{DM02.2})

Here, we focus on the noise effects on blow-up, but in the
non-conservative case, i.e.,
\begin{align} \label{non-conserv}
    \exists j_0: 1\leq j_0 \leq N,\ such\ that\ Re\mu_{j_0} \not =
    0
\end{align}
(Without loss of generality, we assume that $Re\mu_{1} \not =0$).
Surprisingly, the noise effects here are completely different from
those in the conservative cases. We will prove that, in the
non-conservative case by adding a large noise, with high probability
one can prevent blow-up on any given bounded time interval $[0,T]$,
$0<T<\9$. Moreover, when the noise is spatially independent, the
explosion even can be prevented with high
probability on the whole time interval $[0,\9)$. \\

To state our resutls precisely, we assume for the spatial functions
in the noise that
\begin{itemize}
\item[{\rm(H)}] $e_j=f_j+c_j$, $1\leq j \leq N$, where $c_j$
are real constants and $f_j$ are real-valued functions, such that
$f_j\in C_b^{\9}$ and
\begin{align*}
    \lim\limits_{|\xi|\to \9} \zeta(\xi) \sum\limits_{1\leq |\gamma|\leq 3}  |
\partial^{\gamma} f_j(\xi)|  =0,
\end{align*}
where $\gamma$ is a multi-index and
\begin{align*}
    \zeta=\left\{
            \begin{array}{ll}
              1+|\xi|^2, & \hbox{if $d \not =2$;} \\
              (1+|\xi|^2) (\ln (1+|\xi|^2))^2, & \hbox{if $d=2$.}
            \end{array}
          \right.
\end{align*}
(In Section \ref{Main-NEBU} we will take $c_1$ large enough such
that $c_1>|f_1|_{\9}$. Hence, without loss of generality, we assume
that $f_1$ is positive.)
\end{itemize}

The main result is then as follows:
\begin{theorem} \label{Thm}
Consider \eqref{equax} in the non-conservative case
\eqref{non-conserv}. Let $\lbb$ and $\a$ satisfy \eqref{f-super}.
Assume $(H)$ with $f_j$, $1\leq j\leq N$, and $c_k$, $2\leq k\leq N$
being fixed. Then for any $x\in H^1$ and $0<T<\9$,
\begin{align*}
   \mathbb{P} (X(t)\ does\ not\ blow\ up\ \ on\ [0,T]) \to 1,\ \ as\ c_1\to \9.
\end{align*}
(where we recall that by \eqref{non-conserv} we have $Re \mu_1 \not
=0$.)

Furthermore, if $f_j$, $1\leq j\leq N$, are also constants, then for
any $x\in H^1$,
\begin{align*}
   \mathbb{P} (X(t)\ does\ not\ blow\ up\ \ on\ [0,\9)) \to 1,\ \ as\ c_1\to \9.
\end{align*}
\end{theorem}

\begin{remark}
Theorem \ref{Thm} can be viewed as a complement to \cite{BD05}. It
was proved there that in the conservative supercritical case, i.e.,
$Re \mu_j=0$, $1\leq j\leq N$, $\a\in(1+\frac{4}{d}, \9)$ if $d=1,2$
and $\a\in(\frac{7}{3}, 5)$ if $d=3$, the non-degenerate
multiplicative noise can accelerate blow-up with positive
probability (see Theorem $5.1$ in \cite{BD05}). In contrast to
\cite{BD05}, Theorem \ref{Thm} reveals that in the non-conservative
supercritical and also critical cases specified in \eqref{f-super}
with $d\geq 1$, the large multiplicative noise has the effect to
stabilize the system.
\end{remark}

Similar phenomena happen for the deterministic damped nonlinear
Schr\"{o}dinger equation,
\begin{align} \label{damp-NLS}
     i\partial_t u + \Delta u + |u|^{\a-1}u + iau=0,\ \ a>0.
\end{align}
Note that, this equation is analogous to (\ref{weaky-NEBU}) below in
the special case where the noise $W(t)$ is spatially independent and
$\mu_k\in\mathbb{R}$, $1\leq k\leq N$, i.e.
\begin{align*}
    i\partial_t y -\Delta y -e^{(\a-1)ReW(t)}|y|^{\a-1}y + i\wh\mu y
    =0,\ \ \wh\mu>0.
\end{align*}
This similarity indeed indicates the dissipative effects produced by
the multiplicative noise in the non-conservative case.

The global well-posedness of (\ref{damp-NLS}) is proved in
\cite[Theorem $1$]{OT09} (see also \cite[p.98]{SS99}), provided $a$
is large enough, and the proof is based on the decay estimate of
$e^{it\Delta}$ (see \cite[Lemma $4$]{OT09}).

However, since the decay estimates do not necessarily hold for the
general Schr\"odinger-type operator $A(t)$ in \eqref{weaky-NEBU-A},
we employ here quite different arguments based on the contraction
mapping arguments as in \cite{BRZ13, BRZ14}, involving a second
transformation (see \eqref{tra-2} below) and the Strichartz
estimates established in \cite{MMT08}. The advantage of this proof
is that it is also applicable to the case of spatially dependent noise.\\

This article is structured as follows. In Section \ref{Prel-NEBU} we
apply two transformations to reduce the original stochastic equation
\eqref{equax} to a random equation \eqref{equaz-NEBU} below, which
reveals the dissipative effect produced by the noise in the
non-conservative case. Then the non-explosion results in Theorem
\ref{Thm} are established in Section \ref{Main-NEBU}. Furthermore,
we also show that these results do not generally hold with
probability $1$. Finally, the Appendix contains It\^o-formulas for
the Hamiltonian, variance and momentum that are used in the proof.

\section{Preliminaries.} \label{Prel-NEBU}

Following \cite{BRZ13} and \cite{BRZ14}, we apply the rescaling
transformation
\begin{align}  \label{rescal}
X=e^W y
\end{align}
to \eqref{equax} and obtain the random equation
\begin{align}\label{weaky-NEBU}
   &\frac{\partial y}{\partial t}(t, \xi) = A(t)y(t, \xi) - i
    e^{(\alpha-1)Re W(t, \xi)}|y(t, \xi)|^{\alpha-1}y(t, \xi),\\
   & y(0)=x,  \nonumber
\end{align}
where
\begin{equation} \label{weaky-NEBU-A}
A(t) = -i (\Delta  + b(t)\cdot \nabla  + c(t) ),
\end{equation}
\begin{align}
b(t)=2 \nabla W(t),
\end{align}
\begin{align}
c(t)= \sum\limits_{j=1}^d (\partial_j W(t))^2 + \Delta W(t)- i
\widehat{\mu},
\end{align}
and
\begin{align}
   \widehat{\mu}:=\sum\limits_{j=1}^N (|\mu_j|^2 + \mu_j^2)e_j^2.
\end{align}
We stress that on a heuristic level \eqref{weaky-NEBU} follows
easily by It\^o's product rule. The rigorous proof is more involved.
We refer to \cite[Lemma $6.1$]{BRZ13} for the $L^2$-case and
\cite[Theorem $2.1.3$]{Z14} for the $H^1$-case.

Note that the real part of the damped term $\widehat{\mu}$ is
positive in the non-conservative case, namely,
\begin{align} \label{mu-posit}
Re\widehat{\mu}=\sum\limits_{j=1}^{N}(Re\mu_j)^2e_j^2\geq
(Re\mu_{1})^2 c^2_{1}>0,
\end{align}
but it vanishes in the conservative case, which indicates
the different noise effects between the two cases.\\

To explore this damped term, we apply to \eqref{weaky-NEBU} a second
transformation
\begin{align} \label{tra-2}
 z(t,\xi)=e^{\widehat{\mu}
t} y(t,\xi),
\end{align}
and derive that
\begin{align} \label{equaz-NEBU}
   & \frac{\partial z(t)}{\partial t} = \widehat A(t) z(t) -i e^{-(\a-1)(Re\widehat{\mu} t - Re W(t))}  |z(t)|^{\a-1}z(t),\\
   & z(0)=x \in H^1, \nonumber
\end{align}
where \begin{align} \label{Op-A-NEBU}
 \widehat A(t)=-i(\Delta+ \widehat{b}(t)\cdot \na +\widehat {c}(t))
\end{align} with
\begin{align} \label{A-b-NEBU}
\widehat {b}(t)= -2t \na \widehat{\mu} +2\na W(t),
\end{align} and
\begin{align} \label{A-c-NEBU}
      \widehat {c}(t)
     =& t^2 \sum\limits_{j=1}^N(\partial_j
\widehat{\mu})^2 -t\Delta\widehat{\mu}-2t\na W(t) \cdot \na
\widehat{\mu} \nonumber \\
     &+\left[\sum\limits_{j=1}^N (\partial_j W(t))^2 +
\Delta W(t) \right].
\end{align}

The key fact here is that, an exponential decay term $e^{-(\a-1)
Re\widehat{\mu} t} $ appears in \eqref{equaz-NEBU}, which weakens
the nonlinearity and thus can be expected to prevent blow-up,
provided that $\mu$ is sufficiently large (or the noise is
sufficiently large in some other appropriate sense). For this
purpose, let us rewrite equation (\ref{equaz-NEBU}) in the mild form
\begin{equation} \label{mildz-NEBU}
   z(t)=V(t,0)x
   + \int_0^t (-i) V(t,s) \left[ h(s) |z(s)|^{\a-1}z(s) \right]
ds,
\end{equation}
where
\begin{equation} \label{h}
   h(s):=e^{-(\a-1)(Re\widehat{\mu} s - Re W(s))}
\end{equation}
and $V(t,s)$ is the evolution operator generated by the homogenous
part of \eqref{weaky-NEBU}, namely, $V(t,s)x=z(t)$, $s\leq t\leq T$,
solves
\begin{align} \label{equa-v}
   &\frac{d z(t)}{d t} = \widehat A(t) z(t) ,\ \ a.e\ t\in(s,T),\\
   & z(s)=x \in H^1. \nonumber
\end{align}
(The existence and uniqueness of the evolution operator $V(t,s)$
follow mainly from \cite{D94, D96}. For more details, we refer to
\cite{BRZ13,BRZ14}.)

\begin{remark} \label{remark-z}
The solutions to \eqref{equaz-NEBU} are understood analogously to
Definition \ref{def-x}, and Assumption $(H)$ is sufficient to
establish the local existence and uniqueness of solutions for
\eqref{equaz-NEBU}, hence also for \eqref{equax}, by the
transformations \eqref{rescal} and \eqref{tra-2}. Indeed, the proofs
follow by similar arguments as in \cite[Proposition $2.5$]{BRZ14}
(see also \cite[Lemma $4.2$]{BRZ13}), and one can remove the
additional decay assumption $\lim\limits_{|\xi|\to 0} \zeta (\xi)
|e_j(\xi)| =0$ in \cite{BRZ14}, due to the fact that $\wh{b},\wh{c}$
in \eqref{equaz-NEBU} only involve the\ gradient of $\wh \mu$ and
$W(t)$. This fact allows us later to take $c_1$ very large to
prevent blow-up.
\end{remark}

As in \cite[Lemma $2.7$]{BRZ14} one can check from \cite{MMT08} and
Assumption $(H)$ that Strichartz estimates hold for $V(t,s)$,

\begin{lemma} \label{Stri-S}
Assume $(H)$. Then for any $T>0$, $u_0\in H^1$ and $f\in
L^{q_2'}(0,T;W^{1,p_2'})$, the solution of
\begin{equation} \label{stri2}
    u(t)=V(t,0)u_0 + \int_0^t V(t,s) f(s) ds, 0\leq t \leq T,
\end{equation}
satisfies the estimates
\begin{equation} \label{stri-l}
    \|u\|_{L^{q_1}(0,T;L^{p_1})}\leq
    C_T(|u_0|_{2}+\|f\|_{L^{q_2'}(0,T;L^{p_2'})}),
\end{equation}
and
\begin{equation} \label{stri-s}
    \|u\|_{L^{q_1}(0,T;W^{1,p_1})}\leq
    C_T(|u_0|_{H^1}+\|f\|_{L^{q_2'}(0,T;W^{1,p_2'})}),
\end{equation}
where $(p_1,q_1)$ and $(p_2,q_2)$ are Strichartz pairs, i.e.,
\begin{equation*}
    (p_i,q_i)\in[2,\infty] \times
    [2,\infty]: \frac{2}{q_i}=\frac{d}{2}-\frac{d}{p_i},~if~d
    \neq 2,
\end{equation*}
or
\begin{equation*}
    (p_i,q_i)\in[2,\infty) \times
    (2,\infty]:\frac{2}{q_i}=\frac{d}{2}-\frac{d}{p_i},~if~d
    =2,
\end{equation*}
Furthermore, the process $C_t$, $t\geq 0$, can be taken to be
$(\mathscr{F}_t)$-progressively measurable, increasing and
continuous.
\end{lemma}

\section{Proof of the main results} \label{Main-NEBU}

{\it Proof of Theorem \ref{Thm}.} $(i)$. For convenience, let us
first consider the easier case of spatially independent noise to
illustrate the main idea.

By the transformations \eqref{rescal} and \eqref{tra-2}, it is
equivalent to prove the assertion for the random equation
\eqref{equaz-NEBU}. Note that in this case $\wh b=\wh c=0$, hence
$V(t,s)=e^{-i(t-s)\Delta}$ and the Strichartz coefficient $C_t\equiv
C$ is independent of $t$.

Choose the Strichartz pair $(p,q)=(\a+1,\frac{4(\a+1)}{d(\a-1)})$.
Set
\begin{align} \label{damp-NLS.1}
   &\mathcal{Z}^{\tau}_{M}=\{u\in
C(0,\tau;L^2)\cap L^q(0,\tau;L^p):\|u\|_{L^{\infty}(0,\tau; H^1)} +
\|u\|_{L^q(0,\tau; W^{1,p})} \leq M \},
\end{align}
and define the integral operator $G$ on $\mathcal{Z}^{\tau}_{M}$ by
\begin{align} \label{damp-NLS.2}
G(u)(t)=V( t,0) x
    +\int_0^t (-i) V( t, s) \left[h(s) |u(s)|^{\a-1}u(s)\right]
ds,\ u\in \mathcal{Z}^{\tau}_{M}.
\end{align}

We claim that, for $u\in\mathcal{Z}^{\tau}_{M}$,
\begin{align} \label{esti-3.1}
     \|G(u)\|_{L^{\infty}(0,\tau; H^1)} + \|G(u)\|_{L^q(0,\tau; W^{1,p})}
     \leq  2C |x|_{H^1}
    + 2 C D_1(\tau) M^{\a} ,
\end{align}
where
 \begin{align} \label{d3-NEBU}
 D_1(t)= \a D^{\a-1} \|h\|_{L^{v}(0,t)}
\end{align}
with $D$ the Sobolev coefficient such that $\|u\|_{L^p}\leq D
|u|_{H^1}$, $v>1$ and $\frac{1}{v}=1-\frac{2}{q}>0$.

Indeed, by Lemma \ref{Stri-S},
\begin{align} \label{esti-3.1*}
     &\|G(u)\|_{L^{\infty}(0,\tau; H^1)} + \|G(u)\|_{L^q(0,\tau; W^{1,p})}  \nonumber \\
     \leq& 2C |x|_{H^1}
     +2C \|h |u|^{\a-1}u\|_{L^{q'}(0,\tau;W^{1,p'})}.
\end{align}
Moreover, H\"older's inequality and Sobolev's imbedding theorem
yield
\begin{align} \label{esti-3.1.1}
   \|h|u|^{\a-1}u\|_{L^{q'}(0,\tau;L^{p'})}
   \leq& |h|_{L^v(0,\tau)}
   \||u|^{\a-1}u\|_{L^{q}(0,\tau;L^{p'})} \nonumber \\
   \leq& D^{\a-1} |h|_{L^v(0,\tau)}
   \|u\|^{\a-1}_{L^{\9}(0,\tau;H^1)} \|u\|_{L^q(0,\tau;L^{p})},
\end{align}
and
\begin{align} \label{esti-3.1.2}
  \|h\na(|u|^{\a-1}u)\|_{L^{q'}(0,\tau;L^{p'})}
  \leq& \a \|h|u|^{\a-1} |\na u|\|_{L^{q'}(0,\tau;L^{p'})} \nonumber \\
  \leq& \a D^{\a-1} |h|_{L^v(0,\tau)}
   \|u\|^{\a-1}_{L^{\9}(0,\tau;H^1)} \|\na
   u\|_{L^q(0,\tau;L^{p})}.
\end{align}
Hence, plugging \eqref{esti-3.1.1} and \eqref{esti-3.1.2} into
\eqref{esti-3.1*} implies \eqref{esti-3.1}, as claimed. \\

Similarly to \eqref{esti-3.1}, for $u_1,u_2\in
\mathcal{Z}^{\tau}_{M}$,
\begin{align} \label{esti-3.2}
     &\|G(u_1)-G(u_2)\|_{L^{\9}(0,\tau; L^{2})} +  \|G(u_1)-G(u_2)\|_{L^q(0,\tau; L^{p})} \nonumber \\
     \leq& 4 C D_1(\tau) M^{\a-1} \|u_1-u_2\|_{L^q(0,\tau; L^{p})}.
\end{align}

Now, let $M=3C|x|_{H^1}$, choose the $(\mathscr{F}_t)$-stopping time
$\tau=\tau(c_1)$,
\begin{align} \label{cons-tau}
    \tau:=\inf\left\{ t>0:  2 \cdot 3^{\a} |x|_{H^1}^{\a-1} C^{\a} D_1(t) >1
    \right\}.
\end{align}
Then, as in the proof of Proposition $2.5$ in \cite{BRZ14}, we
obtain a local
solution $z$ of \eqref{equaz-NEBU} on $[0,\tau]$.\\

Next we show that $\mathbb{P}(\tau=\9) \to 1$, as $c_1 \to \9$. As
the definition of $\tau$ involves the term $D_1(t)$, we shall use
\eqref{d3-NEBU} to estimate $\|h\|_{L^{v}(0,\9)}$.

Set $\phi_k=\mu_ke_k$, $1\leq k\leq N$. By the scaling property of
Brownian motion, i.e. $\mathbb{P}\circ [Re \phi_k\
\beta_k(\cdot)]^{-1} = \mathbb{P}\circ
[\beta_k((Re\phi_k)^2\cdot)]^{-1}$, for any $c \geq 0$,
\begin{align} \label{cons-hesti}
   &\mathbb{P}(\|h\|^v_{L^{v}(0,\9)} \geq c)  \nonumber \\
   =& \mathbb{P} \( \int_0^{\9} \prod\limits_{k=1}^N e^{-(\a-1)v[(Re\phi_k)^2s-Re\phi_k\beta_k(s)]} ds \geq c
   \)  \nonumber \\
   =&\mathbb{P} \( \int_0^{\9} \prod\limits_{k=1}^N e^{-(\a-1)v[(Re\phi_k)^2s-\beta_k((Re\phi_k)^2 s)]} ds \geq
   c\).
\end{align}

Note that, by the law of the iterated logarithm of Brownian motion,
\begin{align}  \label{int-bdd.2}
 C_1^*:= \int_0^{\9} e^{-(\a-1)v[s-\beta_1(s)]} ds <\9,\ a.s,
\end{align}
and
\begin{align} \label{int-bdd.1}
  C:=1\vee \max\limits_{2\leq k\leq N}\sup\limits_{s\geq
0}e^{-(\a-1)v[(Re\phi_k)^2s-\beta_k((Re\phi_k)^2s)]}< \9,\ a.s.
\end{align}

Then $\mathbb{P}$-a.s.,
\begin{align} \label{cons-intesti}
   &\int_0^{\9} \prod\limits_{k=1}^N e^{-(\a-1)v[(Re\phi_k)^2s-\beta_k((Re\phi_k)^2 s)]}
   ds \nonumber \\
   \leq&   C^N
   \int_0^{\9} e^{-(\a-1)v[(Re\phi_1)^2s-\beta_1((Re\phi_1)^2 s)]} ds  \nonumber \\
   \leq&  \frac{1}{(Re\phi_1)^2}   C^N    C^*_1.
\end{align}
Hence, plugging \eqref{cons-intesti} into \eqref{cons-hesti}, since
$  C^N   C^*_1 <\9\ a.s.$ and $(Re\phi_1)^2 \to \9$ as $c_1 \to \9$,
we deduce that for any fixed $c\geq 0$,
\begin{align} \label{cons-to0}
    &\mathbb{P}(\|h\|^v_{L^{v}(0,\9)} \geq c)
   \leq  \mathbb{P} \(  C^N \wt C^*_1 \geq  c (Re\phi_1)^2\)
    \to 0,\ \ as\ c_1\to \9.
\end{align}

Consequently, choose $c= \left[4 \cdot 3^{\a} \a |x|_{H^1}^{\a-1}
C^{\a} D^{\a-1} \right]^{-v} >0 $. By the definition of $\tau$ in
\eqref{cons-tau} and \eqref{cons-to0}, we then derive that
\begin{align*}
    & \mathbb{P} (\tau=\9)\\
    =& \mathbb{P} \( 2 \cdot 3^{\a} |x|_{H^1}^{\a-1} C^{\a} D_1(t)
    <1,\ \forall t\in[0,\9)\) \\
    \geq& \mathbb{P} \( 2 \cdot 3^{\a} \a |x|_{H^1}^{\a-1} C^{\a}
    D^{\a-1} \|h\|_{L^v(0,\9)} \leq \frac{1}{2} \)\\
    \geq& 1- \mathbb{P} \( \|h\|^v_{L^v(0,\9)} \geq c\)\\
    \to& 1,\ \ as\ c_1\to \9,
\end{align*}
which completes the proof for spatially independent noise.\\

$(ii)$. Now, we consider the general case when the noise $W(t)$ is
space-dependent. Again it is equivalent to prove the assertion for
the random equation \eqref{equaz-NEBU}.

Let $\mathcal{Z}_{M}^{\tau} $, $G$ be as in \eqref{damp-NLS.1} and
\eqref{damp-NLS.2} respectively. Similarly to \eqref{esti-3.1}, for
$u\in \mathcal{Z}_{M}^{\tau} $,
\begin{align} \label{esti-G1}
  & \|G(u)\|_{L^{\infty}(0,\tau; H^1)} +\|G(u)\|_{L^q(0,\tau; W^{1,p})} \nonumber\\
  \leq & 2C_{\tau} |x|_{H^1} +  2C_{\tau} D_2(\tau) M^{\a},
\end{align}
where $C_t$ is the Strichartz coefficient, and
\begin{align}  \label{D}
    D_2(t)
    = \a D^{\a-1}  \|h\|_{L^{v}(0,t;W^{1,\9})}.
\end{align}
with $v>1$ and $\frac{1}{v} =1- \frac{2}{q} >0$.

Moreover, for $u_1,u_2 \in \mathcal{Z}_{M}^{\tau}$,
\begin{align} \label{esti-1.2}
   &\|G(u_1)-G(u_2)\|_{L^{\infty}(0,\tau;L^2)} + \|G(u_1)-G(u_2)\|_{L^q(0,\tau;
   L^{p})} \nonumber \\
   \leq& 4C_{\tau}D_2(\tau) M^{\a-1} \|u_1-u_2\|_{L^q(0,\tau;
    L^{p})}.
\end{align}

Set $M=3C_{\tau} |x|_{H^1}$, choose the $(\mathscr{F}_t)$-stopping
time $\tau=\tau(c_1)$,
\begin{equation} \label{tau1}
   \tau:=\inf \{ t\in[0,T], 2\cdot 3^{\a}|x|^{\a-1}_{H^1} C_t^{\a} D_2(t)  >
1\}\wedge T.
\end{equation}

It follows from \eqref{esti-G1} and (\ref{esti-1.2}) that
$G(\mathcal{Z}_{M}^{\tau}) \subset \mathcal{Z}_{M}^{\tau}$ and $G$
is a contraction on $C([0,\tau];L^2)\cap L^q(0,\tau;L^{p})$.
Therefore, using the same arguments as in \cite{BRZ14}, we obtain a
local
solution $z$ on $[0,\tau]$.\\

To show that $\mathbb{P}(\tau=T) \to 1$, as $c_1 \to \9$, using
\eqref{tau1} and \eqref{D}, we shall estimate $
\|h\|_{L^{v}(0,t;W^{1,\9})}$ below. For simplicity, set $|f|_{\9} :=
|f|_{L^{\9}}$ for any $f\in L^{\9}(\mathbb{R}^d)$ and
$\phi_k:=\mu_ke_k$, $1\leq k \leq N$.

As regards the norm $\|h\|_{L^v(0,t;L^{\9})}$, by \eqref{h} and
\eqref{mu-posit},
\begin{align} \label{tau1.0}
   |h(t)|_{L^{\9}}
   \leq& e^{-(\a-1)\sum\limits_{k=1}^N [\frac{(Re\mu_1)^2c_1^2}{N}t-|Re\phi_k|_{\9}| \beta_k(t)|]
   }.
\end{align}
Analogously to \eqref{int-bdd.1},
\begin{align} \label{tau1.1}
  \wt C:=1\vee \max\limits_{2\leq k\leq N}\sup\limits_{t\geq0} e^{-(\a-1)v[\frac{(Re\mu_1)^2c_1^2}{N}t-|\beta_k(|Re\phi_k|^2_{\9} t)|]}
    < \9,\ \ a.s.
\end{align}

Moreover, choosing $c_1$ large enough such that $c_1>|f_1|_{{\9}}$,
we have
\begin{align} \label{tau1.2}
  &\int_0^T e^{-(\a-1)v[\frac{(Re\mu_1)^2c_1^2}{N}t-|\beta_1(|Re\phi_1|^2_{\9}t)|]} dt  \nonumber \\
  =&\frac{1}{|Re\phi_1|^2_{\9}}
    \int_0^{|Re\phi_1|^2_{\9} T} e^{-(\a-1)v[\frac{(Re\mu_1)^2c_1^2}{N|Re\phi_1|^2_{\9}}
    t-|\beta_1(t)|]}
    dt    \nonumber  \\
  \leq&  \frac{1}{|Re\phi_1|^2_{\9}} \wt C_1^*,
\end{align}
where $\wt C_1^*:=  \int_0^{\9} e^{-(\a-1)v[\frac{1}{4N}t
-|\beta_1(t)|]} dt <\9$ $\mathbb{P}$-a.s.

Thus, as in \eqref{cons-to0}, it follows from
\eqref{tau1.0}-\eqref{tau1.2} and the scaling property of $\beta_k$,
$1\leq k\leq N$, that for any $c>0$ fixed,
\begin{align} \label{tau1.2.1}
   &\mathbb{P}\(C_T^{\a v}\|h\|^v_{L^v(0,T;L^{\9})} \geq c\)  \nonumber \\
   \leq&\mathbb{P} \(C_T^{\a v} \wt C^N \ \wt C_1^* \geq |Re\phi_1|^2_{\9}
   c\)  \nonumber \\
   \to& 0,\ \ as\ c_1\to \9, \ \ \mathbb{P}-a.s.,
\end{align}
where $C_T$ is the Strichartz coefficient.\\

Similar arguments can also be applied to the norm $\|\na
h\|_{L^v(0,t;L^{\9})}$. Indeed, from \eqref{h} and \eqref{mu-posit},
\begin{align*}
   \na h(t)
   =& h(t) \left[-(\a-1) \sum\limits_{k=1}^N \(2Re\phi_k (Re \mu_k \na f_k)
   t-Re\mu_k\na f_k\beta_k(t)\)\right],
\end{align*}
which implies
\begin{align*}
   |\na h(t)|_{\9}
   \leq  (\a-1) |h(t)|_{\9} \sum\limits_{k=1}^N
    \left(2|Re\phi_k|_{\9} |Re\mu_k\na f_k|_{\9} t + |Re\mu_k\na f_k|_{\9}
    |\beta_k(t)|\right).
\end{align*}
Hence, for any $c>0$ fixed,
\begin{align*}
   &\mathbb{P}(C_T^{\a v}\|\na h\|^v_{L^v(0,T;L^{\9})} \geq c)\\
   \leq& \mathbb{P} \big(C_T^{\a v}\int_0^T (\a-1)^v |h(t)|^v_{L^{\9}}\\
   &\qquad\qquad   \left[ \sum\limits_{k=1}^N
       2|Re\phi_k|_{\9} |Re\mu_k\na f_k|_{\9} t + |Re\mu_k\na f_k|_{\9}|\beta_k(t)|
       \right]^v dt \geq c\big)\\
   \leq& \mathbb{P} \bigg(C_T^{\a v}\int_0^T (\a-1)^v
   \left[\prod\limits_{k=1}^N
    e^{-(\a-1)v[\frac{(Re\mu_1)^2c_1^2}{N}t - |\beta_k(|Re\phi_k|^2_{\9} t)|]}\right] \\
      & \qquad\qquad \left[  \sum\limits_{k=1}^N
       2|Re\phi_k|_{\9} |Re\mu_k\na f_k|_{\9} t + |Re\mu_k\na
       f_k|_{\9}|\beta_k(t)|
       \right]^v dt \geq c \bigg)\\
    \leq&  \mathbb{P} \bigg( C_T^{\a v} \wt C^N
         \frac{1}{|Re\phi_1|^2_{\9}} \int_0^{\9} e^{-(\a-1)v [\frac{1}{4N}t - |\beta_1(t)|]}\\
         &\qquad \left[ \sum\limits_{k=1}^N
       \frac{2|Re\phi_k|_{\9} |Re\mu_k\na f_k|_{\9}}{|Re\phi_1|^2_{\9}} t + |Re\mu_k\na f_k|_{\9} \bigg|\beta_k(\frac{t}{|Re\phi_1|^2_{\9}})\bigg|
       \right]^v dt \geq \frac{c}{(\a-1)^v} \bigg).
\end{align*}
Choosing $c_1$ large enough, such that $\sum\limits_{k=1}^N
       \frac{2|Re\phi_k|_{\9} |Re\mu_k\na f_k|_{\9}}{|Re\phi_1|^2_{\9}} < 1$ and
$\frac{|Re\mu_k\na f_k|_{\9}}{|Re\phi_1|_{\9}}<1$, we have as
$c_1\to \9$,
\begin{align} \label{tau1.2.2}
   &\mathbb{P}(C_T^{\a v}\|\na h\|^v_{L^v(0,T;L^{\9})} \geq c) \nonumber \\
   \leq& \mathbb{P}(C_T^{\a v}\wt C^N  \wt C'_1
   \geq \frac{c}{(\a-1)^v}|Re\phi_1|^2_{\9}) \nonumber \\
   \to& 0.
\end{align}
where $C_T$ is the Stichartz coefficient and $ \wt C'_1:=
\int_0^{\9}
e^{-(\a-1)v[\frac{1}{4N}t-|\beta_1(t)|]}\left[t+\sum\limits_{k=1}^N\beta_k(t)\right]^v
dt<\9$ $\mathbb{P}$-a.s. \\

Now we come back to the definition of $\tau$ in \eqref{tau1}.
Choosing $$c=[4\cdot 3^{\a} \a D^{\a-1} |x|_{H^1}^{\a-1}]^{-v}>0,$$
we deduce from \eqref{tau1.2.1} and \eqref{tau1.2.2} that
\begin{align*}
   &\mathbb{P}(\tau=T)\\
   \geq& \mathbb{P}(2\cdot 3^{\a} |x|_{H^1}^{\a-1} C_t^{\a} D_1(t) < 1, \forall t\in[0,T] )\\
   \geq& \mathbb{P}(2\cdot 3^{\a} \a D^{\a-1} |x|_{H^1}^{\a-1} C_T^{\a} \|h\|_{L^v(0,T,W^{1,\9})} < \frac{1}{2})\\
   \geq& 1 -\mathbb{P}(C_T^{\a v} \|h\|^v_{L^v(0,T,W^{1,\9})} \geq c)\\
   \geq&1- \mathbb{P}(C_T^{\a v} \|h\|^v_{L^v(0,T,L^{\9})} \geq \frac{1}{2}c)
       -  \mathbb{P}(C_T^{\a v} \|\na h\|^v_{L^v(0,T,L^{\9})} \geq \frac{1}{2}c)\\
   \to&1, \ \ as\ c_1\to \9.
\end{align*}
Therefore, we complete the proof of Theorem \ref{Thm}. \hfill $\square$\\

One may further ask whether the non-explosion results in Theorem
\ref{Thm} hold with probability $1$. This is, unfortunately, not
generally true. In fact, define the Hamiltonian
\begin{align*}
     H(z) = \frac{1}{2} |\na z|_2^2 - \frac{1}{\a+1}
     |z|_{\a+1}^{\a+1},\ \ z \in H^1,
\end{align*}
and set $\sum=\{u\in H^1, \int |\xi|^2 |u(\xi)|^2 d\xi <\9.\}$. We
have the following result

\begin{proposition}\label{thm2-NEBU}
Consider \eqref{equax} in the non-conservative case
\eqref{non-conserv}. Let $\lbb$ and $\a$ satisfy \eqref{f-super}.
Assume $(H)$ with $f_j$, $1\leq j\leq N$, and $c_k$, $2\leq k\leq N$
being fixed. Furthermore, assume $\mu_k\in\mathbb{R}$, $1\leq k \leq
N$. Let $x\in \sum$ with $H(x)<0$,

Then there exists $\epsilon_0>0$, such that for
$0<\epsilon<\epsilon_0$ and $0\leq \sum\limits_{1\leq k\leq N} |\na
f_k|_{L^{\9}}< \epsilon$, the solution to (\ref{equax}) blows up in
finite time with positive probability.

In particular, in the case that $f_j$, $1\leq j\leq N$, are fixed
constants, the solution to (\ref{equax}) blows up in finite time
with positive probability.
\end{proposition}

The proof follows from the standard virial analysis (see e.g
\cite{G77}). For any $u\in \sum$, define the variance
\begin{align} \label{V}
   V(u)=\int |\xi|^2 |u(\xi)|^2 d\xi,
\end{align}
and the momentum
\begin{align} \label{G}
   G(u)=Im \int \xi u(\xi) \cdot \overline{\na u(\xi)} d\xi.
\end{align}

{\it \bf Proof of Proposition \ref{thm2-NEBU}}. We prove the
assertion by contradiction. Assume that the solution
$X(t)$ to (\ref{equax}) exists globally in $H^1$ $\mathbb{P}-a.s$.\\

By Lemmas $4.1$, \ref{Ito-V} and \ref{Ito-G} in the Appendix,
\begin{align} \label{Ito-V-thm2}
   V(X(t))=&V(x)+4G(x)t+8H(x)t^2 \nonumber \\
   &+4 \sum\limits_{k=1}^N  \int_0^t (t-s)^2 |\na\phi_kX(s)|_2^2ds\nonumber  \\
   &-  4(\a-1) \sum\limits_{k=1}^N  \int_0^t (t-s)^2 \int \phi_k^2|X(s)|^{\a+1}d\xi  ds \nonumber \\
   &+\frac{16}{\a+1}\left[1-\frac{d(\a-1)}{4}\right]\int_0^t(t-s)|X(s)|_{\a+1}^{\a+1}
   ds \\
   &+ M_t \nonumber ,
\end{align}
where $\phi_k=\mu_k e_k$, $1\leq k\leq N$, and
\begin{align*}
   M_t:=&8 \sum\limits_{k=1}^N \int_0^t (t-s)^2 \left[Re\<\na(\phi_kX(s)),\na X(s)\>_2-\int
   \phi_k|X(s)|^{\a+1}d\xi\right] d\beta_k(s)\\
   &-8 \sum\limits_{k=1}^N  \int_0^t (t-s)  Im\int \xi \cdot \na X(s) \overline{X(s)}\phi_k
   d\xi  d\beta_k(s)\\
   &+2 \sum\limits_{k=1}^N  \int_0^t \int |\xi|^2 |X(s)|^2 \phi_k d\xi d\beta_k(s).
\end{align*}

Fix $t>0$ and define for $r\in[0,\9)$,
\begin{align} \label{Ito-V-thm2**}
   \wt M(t,r):=&8 \sum\limits_{k=1}^N \int_0^r (t-s)^2 \left[Re\<\na(\phi_kX(s)),\na X(s)\>_2-\int
   \phi_k|X(s)|^{\a+1}d\xi\right] d\beta_k(s) \nonumber \\
   &-8 \sum\limits_{k=1}^N  \int_0^r (t-s)  Im\int \xi \cdot \na X(s) \overline{X(s)}\phi_k
   d\xi  d\beta_k(s) \nonumber \\
   &+2 \sum\limits_{k=1}^N  \int_0^r \int |\xi|^2 |X(s)|^2 \phi_k d\xi d\beta_k(s).
\end{align}
Set $\sigma_{m}:=\inf\{s\in[0,t],|\na X_m(s)|_2^2>m\}\wedge t$. Then
$\sigma_{m}\to t$, as $m\to\9$.

Direct computations show that $\wt M(t,\cdot \wedge \sigma_{m})$ is
a square integrable martingale, in particular,
\begin{align}
         \mathbb{E} [\wt M(t,t \wedge \sigma_{m})] =0.
\end{align} Indeed, e.g. in regard to the second term in the right hand side
of \eqref{Ito-V-thm2**} , we note that
\begin{align} \label{esti}
    &\mathbb{E} \int_0^{r\wedge \sigma_{m}}  \sum\limits_{k=1}^N \big|(t-s)  Im\int \xi \cdot \na X(s) \overline{X(s)}\phi_k
   d\xi \big|^2  ds \nonumber \\
   \leq& C \mathbb{E} \int_0^{r\wedge \sigma_{m}} (t-s)^2 V(X(s))
   |\na X(s)|_2^2 ds \nonumber \\
   \leq& m C \mathbb{E} \sup\limits_{s\in[0,\sigma_{m}]} V(X(s))
   \int_0^{r} (t-s)^2 ds,
\end{align}
where $C=\sum\limits_{k=1}^N |\phi_k|_{L^{\9}}^2<\9$. Then, as in
the proof of \eqref{bdd-v**} below, we deduce that the right hand
side in \eqref{esti} is finite. The other terms can be estimated
even more easily.

Now, take the expectation in \eqref{Ito-V-thm2}. Since the fifth and
sixth terms in the right hand side of \eqref{Ito-V-thm2} are
non-positive for $\a$ satisfying \eqref{f-super}, it follows that
\begin{align*}
    \mathbb{E} V(X(\sigma_{m} \wedge t ))
   \leq & V(x)+4G(x)(\sigma_{m}\wedge t) +8H(x)(\sigma_{m} \wedge t )^2\\
   &+ 4  \mathbb{E} \int_0^{\sigma_{m} \wedge t } (\sigma_{m} \wedge t -s)^2 \sum\limits_{k=1}^N|\na \phi_k X(s)|_2^2
   ds ,\ \ t<\9.
\end{align*}
Then, taking $m\to\9$, by Fatou's lemma, and since $\na
\phi_k=\mu_k\na f_k$ and $\mathbb{E}|X(t)|_2^2 =|x|_2^2$, we obtain
\begin{align} \label{esti-v}
   \mathbb{E} V(X(t))
   \leq V(x)+4G(x)t +8H(x)t^2
   + at^3
\end{align}
with $$a=\frac{4}{3}\sum\limits_{k=1}^N |\mu_k||\na
f_k|^2_{L^{\9}}|x|_2^2.$$

Let $f(t)$ denote the right hand side of (\ref{esti-v}), i.e.,
\begin{align*}
   f(t) := V(x)+4G(x)t +8H(x)t^2
   + at^3.
\end{align*}

We claim that, if $\sum\limits_{k=1}^N|\na f_k|_{L^{\9}}$ is small
enough, then there exists $T>0$ such that $f(T)<0$. But, taking into
account $ \mathbb{E}V(X(t))\geq 0$ and (\ref{esti-v}), we get a contradiction.\\

It remains to prove the claim. Since
\begin{align*}
   f'(t)= 3at^2 +16H(x)t +4G(x),
\end{align*}
for  $\sum\limits_{k=1}^N|\na f_k|_{L^{\9}}$ small enough, the
discriminant is positive and the largest root of $f(t)$ is
\begin{align} \label{t*}
   t_*:=\frac{2G(x)}{-4H(x)-\sqrt{16(H(x))^2-3aG(x)}}>0.
\end{align}

Note that, proving the claim is equivalent to showing that $f(t_*)<
0$. Since $f'(t^*)=0$, simple computations show that
\begin{align*}
   f(t_*)
   = \frac{8}{3}H(x)t_*^2 +
   \frac{8}{3}G(x)t_* +V(x).
\end{align*}

Since the largest roof of
\begin{align*}
   g(t): =\frac{8}{3}H(x)t^2 + \frac{8}{3}G(x)t +V(x)
\end{align*}
is
\begin{align*}
   \widetilde{t}_*:=\frac{-G(x)-\sqrt{(G(x))^2-\frac{3}{2}H(x)V(x)}}{2H(x)},
\end{align*}
which is independent of $a$. But by \eqref{t*}, $t_*\to \9$  as $a
\to 0$, yielding that $\wt t_* < t_*$ for $a$ small enough, thereby
implying $f(t_*)<0$ and completing the proof. \hfill $\square$ \\

\section{Appendix.} \label{App-NEBU}

This appendix contains the It\^o-formulas for the Hamiltonian,
variance and momentum. As mentioned in Remark \ref{remark-z}, one
can obtain a local solution $X$ to \eqref{equax} on $[0,\tau_n]$,
$n\in\mathbb{N}$, where $\tau_n$ are $(\mathscr{F}_t)$-stopping
times, and $X$ satisfies $\mathbb{P}$-a.s. for any Strichartz pair
$(\rho,\gamma)$,
\begin{align} \label{ito-v*}
   X|_{[0,t]}\in C([0,t];H^1) \cap L^\gamma(0,t;W^{1,\rho}),\ \ t<\tau^*(x)
\end{align}
with $\tau^*(x)=\lim\limits_{n\to\9} \tau_n$.

Let us start with the It\^o-formula for the Hamiltonian $H(X(t))$
proved in \cite[Theorem $3.1$]{BRZ14}.
\begin{theorem} \label{Hami-Ito}
Let $\a$ satisfy \eqref{f-super}. Set $\phi_j:=\mu_j e_j$,
$j=1,...,N$. Then $\mathbb{P}$-a.s
\begin{align*}
     &H(X(t))\\
    =&H(x)
       +\int_0^t Re \<-\nabla(\mu X(s)),\nabla X(s) \>_2ds
       + \frac{1}{2} \sum\limits_{j=1}^N \int_0^t |\nabla (X(s)\phi_j)|_2^2ds        \\
    & -\frac{1}{2}\lambda (\alpha-1)\sum\limits_{j=1}^N \int_0^t \int (Re\phi_j)^2  |X(s)|^{\alpha+1} d\xi ds\\
    & +\sum^N_{j=1} \int_0^t Re \<\nabla(\phi_j
    X(s)),\nabla X(s) \>_2d\beta_j(s)\\
    & -\lambda \sum^N_{j=1} \int_0^t \int Re\phi_j |X(s)|^{\alpha+1} d\xi
    d\beta_j(s),\ \ 0\leq t<\tau^*(x).
\end{align*}
\end{theorem}

The following lemma is concerned with the It\^o-formula for the
variance.
\begin{lemma}  \label{Ito-V}
Let $\sum$ be as in Proposition \ref{thm2-NEBU} and $x\in \sum$.
Then $\mathbb{P}$-a.s. for $t<\tau^*(x)$,
\begin{align} \label{ito-v}
   V(X(t)) = V(x) + 4\int_0^t G(X(s))ds + M_1(t),
\end{align}
where $G$ is as in \eqref{G} and
\begin{align*}
   M_1(t): = 2 \sum \limits_{k=1}^N \int_0^t \int |\xi|^2 |X(s)|^2 Re\phi_k d\xi
   d\beta_k(s)
\end{align*}
with $\phi_k:=\mu_ke_k$, $1\leq k \leq N$.
\end{lemma}

{\it Proof.} The proof is similar to that in \cite[Lemma
$5.1$]{BRZ14} (see also \cite{K10}), hence we just give a sketch of
it below.

Set $\varphi^{\e}:=\varphi \ast \phi_{\e}$ for any locally
integrable function $\varphi$ mollified by $ \phi_{\e}$, where $
\phi_{\e}=\e^{-d}\phi(\frac{x}{\e})$ and $\phi \in
C_c^{\9}(\mathbb{R}^d)$ is a real-valued nonnegative function with
unit integral. Set $V_{\eta}(u)= \int e^{-\eta |\xi|^2} |\xi|^2
|u(\xi)|^2 d\xi$ and $V(u)= \int |\xi|^2 |u(\xi)|^2 d\xi$ for  $u\in
\sum$.

By \eqref{equax} it follows that $\mathbb{P}$-a.s. for every
$\xi\in\mathbb{R}^d$, $t<\tau^*(x)$,
\begin{align} \label{ito-lp.1**}
    (X(t))^{\e}(\xi)&= x^{\e}(\xi)
    + \int_0^t \left[-i\Delta (X(s))^{\e}(\xi)- (\mu X(s))^{\e}(\xi) -
  i (g(X(s)))^{\e}(\xi) \right] ds  \nonumber \\
  &\qquad + \sum\limits_{k=1}^N \int_0^t (X(s) \phi_j)^{\e}(\xi) d\beta_j(s),
\end{align}
where $g(X(s)):=|X(s)|^{\a-1}X(s)$. For simplicity, we set
$X^{\e}(t):=(X(t))^{\e}(\xi)$ and correspondingly for the other
arguments.

Applying the product rule yields $\mathbb{P}$-a.s.
\begin{align*}
      |X^{\e}(t)|^2= &|x^{\e}|^2 -2Re\int_0^t \overline{X}^{\e}(s) i \Delta
      X^{\e}(s)  ds-2Re\int_0^t \overline{X}^{\e}(s)  (\mu X(s) )^{\e}
      ds\\
      &-2Re\int_0^t \overline{X}^{\e}(s)  i[g(X(s))]^{\e} ds + \sum\limits_{k=1}^N \int_0^t |(X(s)\phi_k)^{\e}|^2 ds\\
      &+2\sum\limits_{k=1}^N Re\int_0^t
\overline{X}^{\e}(s)  (X(s) \phi_k)^{\e}
      d\beta_k(s),\ \ t<\tau^*(x).
\end{align*}
Then, integration over $\mathbb{R}^d$ with
$e^{-\eta|\xi|^2}|\xi|^2$, interchanging integrals and integrating
by parts, we have $\mathbb{P}$-a.s. for $t<\tau^*(x)$,
\begin{align} \label{app-veta}
    V_{\eta}(X^{\e}(t))
     =& V_{\eta}(x^{\e})
       +4Im \int_0^t \int e^{-\eta|\xi|^2} (1-\eta |\xi|^2)
      X^{\e}(s)\xi \cdot \na \overline{X^{\e}}(s) d\xi ds \nonumber \\
     &-2 Re \int_0^t \int e^{-\eta|\xi|^2}|\xi|^2 \overline{X}^{\e}(s)(\mu X(s)
     )^{\e}d\xi ds  \nonumber \\
     &-2 Re \int_0^t \int e^{-\eta|\xi|^2}|\xi|^2 \overline{X}^{\e}(s)
     i[g(X(s))]^{\e} d\xi ds   \nonumber  \\
     &+\sum\limits_{k=1}^N  \int_0^t \int e^{-\eta|\xi|^2}|\xi|^2
     |(X(s)\phi_k)^{\e}|^2 d\xi ds  \nonumber \\
     &+2 \sum\limits_{k=1}^N Re\int_0^t\int e^{-\eta|\xi|^2}|\xi|^2  \overline{X}^{\e}(s)  (X(s)
     \phi_k)^{\e} d\xi d\beta_k(s).
\end{align}

As $\sup\limits_{\xi\in \mathbb{R}^d}
e^{-\eta|\xi|^2}[|(1-\eta|\xi|^2)\xi|+|\xi|^2] <\9$, one can take
the limit $\e \to 0$ in (\ref{app-veta}), which leads to
\begin{align} \label{ito-veta}
      V_{\eta}(X(t))
      =&V_{\eta}(x) +4Im \int_0^t \int e^{-\eta|\xi|^2} (1-\eta |\xi|^2)
      X(s)\xi \cdot \na \overline{X}(s) d\xi ds \nonumber \\
      &+2\sum\limits_{k=1}^N \int_0^t\int e^{-\eta|\xi|^2}|\xi|^2  |X(s)|^2 Re\phi_k
      d\xi d\beta_k(s),\ \ t<\tau^*(x).
\end{align}

To pass to the limit $\eta\to 0$, we shall prove that
\begin{align} \label{bdd-v}
    \sup\limits_{s\in[0,\tau_n]}V(X(s)) \leq
    \widetilde{C}(n)<\9,\ \ \mathbb{P}-a.s.
\end{align}
Then by \eqref{ito-veta}, \eqref{bdd-v},
$\sup\limits_{\eta>0}\sup\limits_{\xi\in\mathbb{R}^d}
|e^{-\eta|\xi|^2}(1-\eta|\xi|^2)|=1$ and Lebesque's dominated
theorem, we obtain \eqref{ito-v} for $t\leq \tau_n$,
$n\in\mathbb{N}$. Consequently, since $\tau_n\to \tau^*(x)$, as
$n\to \9$, we conclude
\eqref{ito-v} for $t<\tau^*(x)$. \\

It remains to prove (\ref{bdd-v}). For every $n\in\mathbb{N}$, set
$$\sigma_{n,m}:=\inf\{s\in[0,\tau_n]: |\na X(s)|_2^2 >m\}\wedge
\tau_n.$$ Burkholder-Davis-Gundy's inequality implies that
\begin{align} \label{bdd-v.0}
     \mathbb{E} \sup\limits_{s\in[0,t\wedge \sigma_{n,m}]} V_{\eta}(X(s))
     \leq & 4 \mathbb{E} \int_0^{t\wedge \sigma_{n,m}} \int e^{-\eta|\xi|^2} |1-\eta |\xi|^2|
      |\overline{X}(s) \xi \cdot \na X(s)| d\xi ds \nonumber \\
      &+ c \mathbb{E}\sqrt{\int_0^{t\wedge \sigma_{n,m}}\sum\limits_{k=1}^N  \(\int e^{-\eta|\xi|^2}|\xi|^2  |X(s)|^2 Re\phi_k
      d\xi\)^2ds} \nonumber  \\
      =&J_1+J_2,
\end{align}
where $c$ is independent of $n$, $m$ and $\eta$.

Since $\sup\limits_{\eta>0}\sup\limits_{\xi\in\mathbb{R}^d}
|e^{-\eta|\xi|^2}(1-\eta|\xi|^2)|=1$ and  $\mathbb{E}
\sup\limits_{s\in[0,\sigma_{n,m}]}
   |\na X(s)|_2^2 \leq m <\9$,
\begin{align} \label{bdd-v.1}
   J_1 \leq& 4 \mathbb{E} \int_0^{t\wedge \sigma_{n,m}} \sqrt{V(X(s))} |\na X(s)|_2
   ds \nonumber \\
   \leq& 4\int_0^t \mathbb{E} \sup\limits_{r\in[0,s\wedge \sigma_{n,m}]} V(X(r)) ds
   + 4mT.
\end{align}
Moreover,
\begin{align} \label{bdd-v.2}
   J_2 \leq& C \mathbb{E} \sqrt{\int_0^{t\wedge \sigma_{n,m}}
   [V_{\eta}(X(s))]^2 ds} \nonumber \\
   \leq& \e C
   \mathbb{E}\sup\limits_{s\in[0,t\wedge \sigma_{n,m}]}V_{\eta}(X(s))
   +C C_{\e}  \int_0^t \mathbb{E} \sup\limits_{r\in[0,s\wedge \sigma_{n,m}]}
   V_{\eta}(X(r))ds,
\end{align}
where $C$ depends on $|\phi_k|_{L^{\9}}$, $1\leq k\leq N$, and is
independent of $n,m$ and $\eta$.

Hence, plugging \eqref{bdd-v.1} and \eqref{bdd-v.2} into
\eqref{bdd-v.0}, taking $\e$ small enough, and noting that
$V_{\eta}(X)\leq V(X)$, we derive that
\begin{align*}
   \mathbb{E}\sup\limits_{s\in[0,t\wedge \sigma_{n,m}]}V_{\eta}(X(s))
   \leq& C_1 \int_0^{  t} \mathbb{E} \sup\limits_{r\in[0,s\wedge \sigma_{n,m}]} V(X(r)) ds
   +C_2(m,T),
\end{align*}
with $C_1$ and $C_2(m,T)$ independent of $\eta$. Then letting
$\eta\to 0$ and using Fatou's lemma, we have
\begin{align*}
     \mathbb{E}\sup\limits_{s\in[0,t\wedge \sigma_{n,m}]}V(X(s))
     \leq C_1 \int_0^{  t} \mathbb{E} \sup\limits_{r\in[0,s\wedge \sigma_{n,m}]} V(X(r)) ds
   +C_2(m,T),\ \ t\in[0,T],
\end{align*}
which implies by Gronwall's inequality that
\begin{align} \label{bdd-v**}
    \mathbb{E} \sup\limits_{t\in[0,\sigma_{n,m}]}V(X(t)) \leq
    C(m,T)<\9,
\end{align}
hence $\sup\limits_{t\in[0,\sigma_{n,m}]}V(X(t)) \leq
   \wt C(m,T)<\9$, $\mathbb{P}$-a.s. But, since $\sup\limits_{t\in[0,\tau_n]} |\na X(t)|_2^2 < \9$,
$\mathbb{P}$-a.s, for $\mathbb{P}$-a.e. $\omega\in\Omega$, $\exists
m (\omega)<\9$ such that $\sigma_{n,m(\omega)}(\omega)
=\tau_n(\omega)$. Then $\mathbb{P}\(\bigcup\limits_{m\in\mathbb{N}}
\{\sigma_{n,m}=\tau_n\} \) =1$. This implies \eqref{bdd-v} and
completes the proof of Lemma \ref{Ito-V}.  \hfill $\square$\\

We conclude this section with the It\^o-formula for the momentum.
\begin{lemma} \label{Ito-G}
Let $x\in \sum$. Then $\mathbb{P}$-a.s for $t<\tau^*(x)$,
\begin{align} \label{ito-g}
   G(X(t))=& G(x) + 4\int_0^t P(X(s)) ds \nonumber \\
     &- \sum\limits_{k=1}^N  \int_0^t Im \int  \xi \cdot\na\phi_k |X(s)|^2
   \overline{\phi_k}  d\xi ds + M_2(t),
\end{align}
where
\begin{align*}
    P(X):=&\frac{1}{2}|\na X|_2^2 -
\frac{d(\a-1)}{4(\a+1)} |X|_{L^{\a+1}}^{\a+1} \\
        =& H(X) +
\frac{1}{\a+1}[1-\frac{d(\a-1)}{4}]|X|_{L^{\a+1}}^{\a+1},
\end{align*}
$\phi_k=\mu_ke_k$, $1\leq k \leq N$, and
\begin{align*}
   M_2(t):=&d \sum\limits_{k=1}^N \int_0^t\int |X(s)|^2 Im\phi_k d\xi
   d\beta_k(s) \\
   &-2 \sum\limits_{k=1}^N \int_0^t Im \int \xi \cdot \na X(s) \overline{X(s)} \overline{\phi_k}
   d\xi d\beta_k(s).
\end{align*}
Here, $d$ is the dimension of the space.
\end{lemma}

{\it \bf Proof.}  The proof is similar to that in Lemma \ref{Ito-V}
but involves more complicated computations. For simplicity of
exposition, we omit the proof here and refer to
\cite[Lemma $3.3.2$]{Z14} for details. \hfill $\square$ \\

{\it \bf Acknowledge.}  Deng Zhang would like to thank DFG for the
support through the IRTG 1132 at Bielefeld University.

\end{document}